\documentclass[reqno,12pt]{amsart}
\usepackage{amsmath,amsthm,amssymb,amsfonts,amscd}
\newtheorem{thm}{Theorem}[section]
\newtheorem{cor}[thm]{Corollary}
\newtheorem{lem}[thm]{Lemma}
\newtheorem{prop}[thm]{Proposition}
\newtheorem{conj}[thm]{Conjecture}

\theoremstyle{definition}
\newtheorem{defn}[thm]{Definition}

\newtheorem{rem}[thm]{Remark}

\newtheorem*{pf}{Proof}
\numberwithin{equation}{section}
\def\C{{\mathbb C}}

\def\Z{{\mathbb Z}}

\def\CC{{\mathcal C}}

\def\L{{\mathcal L}}

\begin{document}
\title{Weighted Projective Lines Associated to\\
Regular Systems of Weights of Dual Type}
\date{\today}
\author{Atsushi Takahashi}
\address{Department of Mathematics, Graduate School of Science, Osaka University, 
Toyonaka Osaka, 560-0043, Japan}
\email{takahashi@math.sci.osaka-u.ac.jp}
\begin{abstract}
We associate to a regular system of weights 
a weighted projective line over an algebraically closed field 
of characteristic zero in two different ways.
One is defined as a quotient stack via a hypersurface singularity for 
a regular system of weights and the other is defined via the signature of 
the same regular system of weights.
The main result in this paper is that if a regular system of weights is of dual type then 
these two weighted projective lines have equivalent abelian categories of coherent sheaves.
As a corollary, we can show that the triangulated categories of the graded singularity associated 
to a regular system of weights has a full exceptional collection, which is expected from 
homological mirror symmetries. 
Main theorem of this paper will be generalized to more general one, to the case when a regular system of weights is 
of genus zero, which will be given in \cite{kst:3}.
Since we need more detailed study of regular systems of weights and some knowledge of algebraic 
geometry of Deligne--Mumford stacks there, the author write a part of the result in this paper to which 
another simple proof based on the idea by Geigle--Lenzing \cite{gl:2} can be applied.
\end{abstract}
\maketitle
\pagestyle{myheadings}
\markleft{Atsushi TAKAHASHI}
\markright{Weighted Projective Lines Associated to Regular Systems of Weights of Dual Type}
\section{Introduction}
Let $W:=(a_1,a_2 ,a_3;h)$ be a tuple of four positive integers.
If it satisfies a certain combinatorical condition, it is called a regular system of weights.
Let $k$ be an algebraically closed field of characteristic zero.
The condition is equivalent to the condition that 
for a polynomial in $k[x,y,z]$
$$
f_W(x,y,z)=\sum_{a_1i_1+a_2i_2+a_3i_3=h}c_{i_1i_2i_3}x^{i_1}y^{i_2}z^{i_3},\quad c_{i_1i_2i_3}\in k,  
$$
with generic coefficients $c_{i_1i_2i_3}$, 
${\rm Spec}(k[x,y,z]/(f_W))$ has at most an isolated singularity only at the origin \cite{s:1}.
Since $f_W$ is weighted homogeneous, one naturally associates to $W$ the following quotient stack:    
$$
\CC_W:=\left[{\rm Spec}(R_{W})\backslash\{0\}/k^*\right],
$$
where we set $k^*:={\rm Spec}(k\Z)$. 
$\CC_W$ is a Deligne--Munford stack regarded as a smooth projective curve ${\rm Proj}(R_W)$ with 
a finite number of isotropic points on it.
On the other hand, to the signature $A_W=(\alpha_1,\dots, \alpha_r)$ of $W$, 
a combinatorically defined multi-set which can be identified with the multi-set of orders of isotropy group at isotropic points 
on ${\rm Proj}(R_W)$ \cite{s:2}, one can associate the algebra following Geigle--Lenzing \cite{gl:1}
$$
R_{A_W,\lambda}:=k[X_1,\dots,X_r]\left/I_{\lambda}\right.,
$$
where  $I_\lambda$ is an ideal genarated by $r-2$ homogeneous polynomials
$$
X_1^{\alpha_1}+X_2^{\alpha_2}+X_3^{\alpha_3}, X_1^{\alpha_1}+\lambda_iX_2^{\alpha_2}+X_i^{\alpha_i}, \quad
\lambda_i\in k\backslash\{0,1\},\quad i=4,\dots, r.
$$
Since $R_{A_W,\lambda}$ is graded with respect to an abelian group
$$
L(A_W):=\bigoplus_{i=1}^r\Z\vec{X}_i\left/\left(\alpha_i\vec{X}_i-\alpha_j\vec{X}_j;1\le i<j\le r\right)\right. ,
$$
one can consider the quotient stack
$$
\CC_{A_W,\lambda}:=\left[{\rm Spec}(R_{A_W,\lambda})\backslash\{0\}/{\rm Spec}({kL(A_W)})\right].
$$
$\CC_{A_W,\lambda}$ is a Deligne--Mumford stack whose underlying quotient scheme is 
a smooth projective line (it is easy to see that $R_{A_W,\lambda}\supset k[X_1^{\alpha_1},X_2^{\alpha_2}]$ as a subring).
Properties of categories ${\rm coh}(\CC_{A_W})$ and $D^b{\rm coh}(\CC_{A_W})$ are extensively studied by 
Geigle--Lenzing \cite{gl:1}.
It is very natural and intersting problem to compare $\CC_W$ with $\CC_{A_W,\lambda}$ as algebraic stacks.
We have the following result.
\begin{thm}[Main Theorem \ref{thm:main}]
Let $W$ be a regular system of weights of dual type. 
One can choose a weighted homogeneous polynomial $f_W$ so that 
it induces an equivalence of abelian categories$:$
\begin{equation}
{\rm coh}(\CC_W) \simeq  {\rm coh}(\CC_{A_W}).
\end{equation}\qed
\end{thm}
One of motivations of this work is to study a qualitative structure of the triangulated category
$D^{gr}_{Sg}(R_W):=D^b({\rm gr}\text{-}R_W)/D^b({\rm grproj}\text{-}R_W)$. 
In particular, we have been intersted in finding a "nice" or "special" full strongly exceptional collection 
of $D^{gr}_{Sg}(R_W)$ by using graded matrix factorizations.
So far, we have succeeded to obtain it in \cite{t:2} for $A_l$-type, in \cite{kst:1} for $ADE$-type and in 
\cite{kst:2} for the case $\epsilon_W=-1$.
It is in general a very difficult problem, however, by the main theorem in this paper combining with 
a result by Orlov \cite{o:1}, one has a slightly weaker statement for {\it any} regular system of weights of dual type: 
\begin{cor}[Corollary \ref{cor:exceptional}]
Let $W$ be a regular system of weights of dual type.
The triangulated category $D^{gr}_{Sg}(R_W)$ 
has a full exceptional collection $(E_1,\dots,E_{\mu_{W^*}})$.
\qed\end{cor}
This is expected from the homological mirror symmetry conjecture that predicts the existence of 
a triangulated equivalence $D^{gr}_{Sg}(R_W)\simeq D^b{\rm Fuk}(f_{W^*})\simeq D^b{\rm Fuk}^{\rightarrow}(\{\gamma_\bullet \})$
where ${\rm Fuk}^{\rightarrow}(\{\gamma_\bullet\})$ is the directed Fukaya category associated to a distinguished basis of 
vanishing graded Lagrangian submanifolds $\{\gamma_\bullet\}=\{\gamma_1,\dots, \gamma_\mu\}$ in 
the Milnor fiber of the dual regular system of weights (see Section \ref{sec:conj} for details).
\smallskip
We give here an outline of the paper.
Section \ref{sec:rws} introduces the definition of a regular system of weights and 
several invariants of it given in \cite{s:1} and \cite{s:2}. 
Section \ref{sec:duality} explains the notion of the topological mirror symmetry and the duality 
of regular systems of weights defined in \cite{t:1} and \cite{s:2}. 
Most of results in this paper rely on the classification of regular systems of weights of dual type 
and several data given explicitly by them. 
In Section \ref{sec:wpl}, after preparing some defintions and reviewing a construction  
of weighted projective line by Geigle--Lenzing \cite{gl:1}, we characterize a special subring $R_W$
of their homogeneous coordinate ring $R_{A_W}$ of a weighted projective line, which will play a key role in our story.
Section \ref{sec:thm} gives the main theorem of this paper. Its proof uses the data of invariants of 
regular systems of weights of dual type, which we gave in Appendix.
We shall give an application of this result in Section \ref{sec:cor}.
Section \ref{sec:conj} explains the homological mirror symmetry conjecture of hypersurface singularities and 
one of our motivations of the paper.
In particular, it is discussed that the bounded derived category of coherent sheaves on a weighted projective line
associated to a regular system of weights of dual type is expected to be triangulated equivalent to 
the derived category of a directed Fukaya category associated to a cusp singularity associated to the dual signature. 
\smallskip
Main theorem of this paper will be generalized to more general one, to the case when a regular system of weights $W$ 
is of genus zero, i.e., when $g({\rm Proj}(R_W))=0$, which will be given in \cite{kst:3}.
\smallskip
{\bf Acknowledgment}.\  
The author would like to thank Hiroshige Kajiura and Kyoji Saito for valuable discussion. 
This work was partly supported by Grant-in Aid for Scientific Research 
grant numbers 17740036 from the Ministry of Education, Culture, Sports, Science and Technology, 
Japan. 
\section{Regular Systems of Weights}\label{sec:rws}
In this paper, we denote by $k$ an algebraically closed field of characterictic zero.
\begin{defn}[\cite{s:1}]
Let $a_i$, $i=1,2,3$ and $h$ be positive integers.
\begin{enumerate}
\item We call a tuple of integers $W:=(a_1,a_2 ,a_3;h)$ a {\it regular system of weights} 
if $\gcd(a_1,a_2 ,a_3;h)=1$ and a rational function:
\begin{equation}
\chi_W(T):=\prod_{i=1}^{3}\frac{1-T^{h-a_i}}{1-T^{a_i}},
\end{equation}
is a polynomial in $T$.
\item The integers $a_i$ are called {\it weights} of $W$ and $h$ is called
the {\it coxeter number} of  $W$.
\end{enumerate}
\end{defn}
The next proposition relates combinatorics of regular systems of weights with geometries of 
hypersurface singularities.
\begin{prop}[\cite{s:1}]
The following conditions are equivalent$:$
\begin{enumerate}
\item $W=(a_1,a_2 ,a_3;h)$ is a regular system of weights. 
\item There is at least one weighted homogeneous polynomial in $k[x,y,z]$ 
$$
f(x,y,z)=\sum_{a_1i_1+a_2 i_2+a_3 i_3=h}c_{i_1i_2i_3}x^{i_1}y^{i_2}z^{i_3},\quad c_{i_1i_2i_3}\in k,
$$
such that the hypersurface $\{(x,y,z)\in k^3~|~f(x,y,z)=0\}$ has at most an isolated singularity only at the origin.
\item There is a non-empty dence subset of  
$$
\{g(x,y,z)\in k[x,y,z]~|~g(x,y,z)=\sum_{a_1i_1+a_2 i_2+a_3 i_3=h}c_{i_1i_2i_3}x^{i_1}y^{i_2}z^{i_3}, c_{i_1i_2i_3}\in k,\}
$$
such that any polynomial $f(x,y,z)$ belonging to that set defines the hypersurface $\{(x,y,z)\in k^3~|~f(x,y,z)=0\}$ 
which has at most an isolated singularity only at the origin.
\end{enumerate}
\qed
\end{prop}
\begin{rem}
We shall see later that one can somtimes choose in $(iii)$ a "canonical" polynomial, 
which will be one of keys to prove the main theorem in this paper.
\end{rem}
\begin{defn}[\cite{s:1}]
Let $W$ be a regular system of weights.
\begin{enumerate}
\item The positive integer $\mu_W$ defined by 
\begin{equation}
\mu_W:=\prod_{i=1}^3\frac{h-a_i}{a_i}
\end{equation}
is called the {\it rank} or the {\it Milnor number} of the regular system of weights $W$.
\item The integer $\epsilon_W$ defined by 
\begin{equation}
\epsilon_W:=\left(\sum_{i=1}^3a_i\right)-h
\end{equation}
is called the {\it minimal exponent} or the {\it Gorenstein parameter} of $W$.
\item There are finite number of
integers $m_1<m_2\le\dots\le m_{\mu_W-1}<m_{\mu_W}$ such that 
$$
T^{m_1}+\dots+T^{m_{\mu_W}}=T^{\epsilon_W}\chi_W(T).
$$
Each integer $m_i$ is called an {\it exponent} of $W$, which satirfies the property
\begin{equation}
m_i+m_{\mu_W-i+1}=h,\quad i=1,\dots ,\mu_W.
\end{equation}
\end{enumerate}
Define more "geometric" invariants for regular systems of weights, whose meanings will be 
clear in later sections.
\end{defn}
\begin{defn}[\cite{s:1}\cite{s:2}]
Let $W=(a_1,a_2,a_3;h)$ be a regular system of weights.
\begin{enumerate}
\item The {\it genus} of $W$ is a non-negative integer $g_W$ defined as the number of exponents of $W$ 
equal to $0$.
\item Set $m(a_i,a_j:h):=\#\{(u,v)\in(\Z_{\ge 0})^2\, |\, a_iu+a_jv=h \}$.
Consider the following multi-set of positive integers
$$
A'_W:= \{a_i~|~h/a_i \notin\Z,i=1,2,3\} \coprod
\{{\rm gcd}(a_i,a_j)^{\left(m(a_i,a_j:h)-1\right)}~|~1\le i<j\le 3\},
$$
where we mean by $\gcd(a_i,a_j)^{\left(m(a_i,a_j:h)-1\right)}$ that the integer ${\rm gcd}(a_i,a_j)$ appears  
$m(a_i,a_j:h)-1$ times in $A_W'$. 
Let $A_W$ be a subset of $A'_W$ consists of integers greater than $1$ 
such that $A'_W=A_W\coprod \{1^s\}$ for some $s\in\Z_{\ge 0}$.
\item The pair $(g_W;A_W)$ is called the {\it signature} of $W$. 
If $g_W$ is $0$, then we often omit $g_W$ and call $A_W$ the signature of $W$. 
\end{enumerate}
\end{defn}
\begin{rem}
Usually, we shall write $A_W$ as $(\alpha_1,\dots,\alpha_r)$ so that $2\le \alpha_1\le\alpha_2\le\cdots\le\alpha_r$ 
for some $r\in\Z_{\ge 0}$.  
\end{rem}
\begin{defn}[\cite{s:2}]
Let $W$ be a regular system of weights.
\begin{enumerate}
\item The polynomial 
\begin{equation}
\varphi_W(\lambda):=\prod_{i=1}^{\mu_W}(\lambda-{\bf e}[\frac{m_i}{h}]),\quad \quad {\bf e}[*]=\exp(2\pi\sqrt{-1}*),
\end{equation}
is called the {\it characteristic polynomial} of $W$.
\item A set of integers $M(W)$ such that 
\begin{equation}
\varphi_W(\lambda)=\prod_{i\in M(W)}(\lambda^i-1)^{e_W(i)}
\end{equation}
is called the {\it classifying poset} of $W$.
\item The {\it dual characteristic polynomial} $\varphi^*_W(\lambda)$ of $W$ is a polynomial 
\begin{equation}
\varphi^*_W(\lambda):=\prod_{i\in M(W)}(\lambda^i-1)^{-e_{W}(h/i)}.
\end{equation}
\end{enumerate}
\end{defn}
\begin{rem}
$M(W)$ is always well-defined since $\varphi_W(\lambda)$ is a cyclotomic polynomial.
It is a partially ordered set (poset) with respect to the division relation on integers. 
\end{rem}
\section{Duality of Regular Systems of Weights}\label{sec:duality}
In this subsection, we recall the notion of the duality of regular system of weights discussed in \cite{s:2} and \cite{t:1}.
\begin{defn}\label{def:vafa}
Let $W$ be a regular system of weights and $G$ be a finite abelian subgroup of $GL(3,k)$, whose elements are of the form
${\rm diag}({\bf e}[\omega_1\alpha_1],{\bf e}[\omega_2\alpha_2] ,{\bf e}[\omega_3\alpha_3])$ where $\omega_i:=a_i/h$ and
$\alpha_i\in\Z$.
For a pair $(W,G)$, set
\begin{equation}
\chi(W,G)(y,\bar{y}):=\frac{(-1)^{3}}{|G|}\sum_{\alpha\in G}\chi_\alpha(W,G)
(y,\bar{y}),
\end{equation}
\begin{equation}
\begin{split}
\chi_\alpha(W,G)(y,\bar{y}):=& \sum_{\beta\in G}
\prod_{\omega_i\alpha_i\not\in\Z}\left(y\bar{y}\right)^{\frac{1-2\omega_i}{2}}
\left(\frac{y}{\bar{y}}\right)
^{-\omega_i\alpha_i+\left[\omega_i\alpha_i\right]+\frac{1}{2}}\\
& \times\prod_{\omega_i\alpha_i\in\Z}{\bf e}\left[\omega_i\beta_i+\frac{1}{2}
\right]\frac{1-{\bf e}\left[(1-\omega_i\beta_i)\right](y\bar{y})^{1-\omega_i}}
{1-{\bf e}\left[\omega_i\beta_i\right](y\bar{y})^{\omega_i}},
\end{split}
\end{equation}
where $[\omega_i\alpha_i]$ denotes the greatest integer smaller than $\omega_i\alpha_i$.
We call $\chi(W,G)(y,\bar{y})$ the {\it orbifoldized Poincare polynomial} of $(W,G)$.
\end{defn}
\begin{rem}
Note that if we put $T^h=y\bar{y}$, we have
$$
\chi_W(T)=\chi(W,\{1\})(y,\bar{y}).
$$
\end{rem}
By this orbifoldized Poincare polynomial $\chi(W,G)$, one can define the following 
notion of the duality among pairs $(W,G)$:
\begin{defn}[Topological Mirror Symmetry]\label{defn:topological mirror}
Let $(W,G)$ and $(W^*,G^*)$ be pairs as in Definition \ref{def:vafa}.
A pair $(W^*,G^*)$ is called {\it topological mirror dual} to a pair $(W,G)$, if
\begin{equation}
\chi(W^*,G^*)(y,\bar{y})=(-1)^3\bar{y}^{\hat{c}_{W}}
\chi(W,G)(y,\bar{y}^{-1}),
\end{equation}
where $\hat{c}_W:=1-2\frac{\epsilon_W}{h}$.
\end{defn}
Let $W$ be a regular system of weights. We shall call the group generated by
${\rm diag}({\bf e}[\omega_1],{\bf e}[\omega_2] ,{\bf e}[\omega_3])$
the {\it principal discrete group} for $W$ and denote it by $G^0_W$.
As a special case of topological mirror symmetry, we define the following duality of regular systems of weights:
\begin{defn}[\cite{t:1}]\label{defn:dual}
Let $W$ and $W^*$ be regular systems of weights.
$W^*$ is said to be {\it dual} to $W$ if
\begin{equation}
\chi(W^*,\{1\})(y,\bar{y})=(-1)^3\bar{y}^{\hat{c}_{W}}
\chi(W,G_W^0)(y,\bar{y}^{-1}).
\end{equation}
\end{defn}
\begin{prop}[\cite{t:1}]
Let $W$ and $W^*$ be regular systems of weights. 
$W^*$ is dual to $W$ if and only if $W^*$ is $*$-dual to $W$ in the sense of K.~Saito \cite{s:2}, 
more precisely, if and only if 
\begin{equation}
\varphi_{W^*}(\lambda)=
\prod_{i\in M(W^*)}(\lambda^i-1)^{e_{W^*}(i)}
=\prod_{i\in M(W)}(\lambda^i-1)^{-e_{W}(h/i)}=\varphi_W^*(\lambda),
\end{equation}
and $W^*=(lm-l+1,mk-m+1,kl-k+1;klm+1)$ when $W=(lm-m+1,mk-k+1,kl-l+1;klm+1)$.
\qed
\end{prop}
\begin{rem}
If $W^*$ is dual to $W$, then $W^*$ is uniquely determined by $W$ and $(W^*)^*=W$. 
Therefore, one sees that these two equivalent notions of duality define equivalence relations. 
\end{rem}
\begin{defn}
A regular system of weights $W$ is called of {\it dual type} if $W$ has the dual regular system of weights $W^*$. 
\end{defn}
\begin{prop}[\cite{s:2}]
Let $W$ be a regular system of weights $W$ of dual type.
Then, the signature $(g_W;A_W)$ of $W$ is of the form  $(0;\alpha_1,\alpha_2,\alpha_3)$.
\qed
\end{prop}
\begin{rem}
Regular systems of weights of dual type is classified into five types with respect to the classifying poset $M(W)$.
See Appendix for details of their data, i.e., weights, dual weights, signatures, and so on.  
\end{rem}
\section{Weighted Projective Lines Associated to $A_W$}\label{sec:wpl}
\begin{defn}[\cite{gl:1}]\label{defn:L}
Let $W=(a_1,a_2,a_3;h)$ be a regular system of weights of genus $0$ and $A_W=(\alpha_1,\dots, \alpha_r)$ be its signature. 
\begin{enumerate}
\item Consider an abelian group generated by $r$-letters $\vec{X_i}$ $(i=1,\dots ,r)$, 
\begin{equation}
L(A_W):=\bigoplus_{i=1}^r\Z\vec{X}_i\left/\left(\alpha_i\vec{X}_i-\alpha_j\vec{X}_j;1\le i<j\le r\right)\right. .
\end{equation}
It is an ordered group with $L(A_W)_+:=\sum_{i=1}^r\Z_{>0}\vec{X_i}$ as its positive elements.
\item An element $\omega_{A_W}:=(r-2)\cdot\vec{c}-\sum_{i=1}^r\vec{X_i}$ is called the {\it dualizing elment} of $L(A_W)$, 
where $\vec{c}:=\alpha_1\vec{X_1}=\dots=\alpha_r\vec{X_r}$.
\item Let $\alpha:={\rm lcm}(\alpha_1,\dots,\alpha_r)$. 
The {\it degree map} is a group homomorpshim 
${\rm deg}:L(A_W)\longrightarrow \Z$ defined on generators by ${\rm deg}(\vec{X_i}):=\alpha/\alpha_i$.
The degree map ${\rm deg}$ is an epimorphism and ${\rm deg}(\vec{X})=0$ if $\alpha\cdot\vec{X}=0$. 
Note that ${\rm deg}(\vec{c})=\alpha$.
\item Define an $L(A_W)$-graded $k$-algebra by 
\begin{equation}
R_{A_W,\lambda}:=k[X_1,\dots,X_r]\left/I_{\lambda}\right.,
\end{equation}
where $I_\lambda$ is an ideal genarated by $r-2$ homogeneous polynomials
\begin{equation}
X_1^{\alpha_1}+X_2^{\alpha_2}+X_3^{\alpha_3}, X_1^{\alpha_1}+\lambda_iX_2^{\alpha_2}+X_i^{\alpha_i}, \quad
\lambda_i\in k\backslash\{0,1\},\quad i=4,\dots, r.
\end{equation}
$R_{A_W,\lambda}$ is non-negatively graded.
If $r=3$, then we shall write $R_{A_W,\lambda}$ as $R_{A_W}$ for simplicity.
\end{enumerate}
\end{defn}
Denote by ${\rm mod}_{L(A_W)}\text{-}R_{A_W,\lambda}$ the abelian category of finitely generated $L(A_W)$-graded 
$R_{A_W,\lambda}$-modules and denote by ${\rm tor}_{L(A_W)}\text{-}R_{A_W,\lambda}$ the full subcategory of 
${\rm mod}_{L(A_W)}\text{-}R_{A_W,\lambda}$ whose objects are finite dimensional $L(A_W)$-graded $R_{A_W,\lambda}$-modules.
\begin{defn}
Define an algebraic stack $\CC_{A_W,\lambda}$ by 
\begin{equation}
\CC_{A_W,\lambda}:=\left[{\rm Spec}(R_{A_W,\lambda})\backslash\{0\}/{\rm Spec}({k\cdot L(A_W)})\right].
\end{equation}
The abelian category ${\rm coh}(\CC_{A_W,\lambda})$ of coherent sheaves on the algebraic stack 
$\CC_{A_W,\lambda}$ is equivalent to 
${\rm mod}_{L(A_W)}\text{-}R_{A_W,\lambda}/{\rm tor}_{L(A_W)}\text{-}R_{A_W,\lambda}$.
\end{defn}
$L(A_W)$ has a special subgroup isomorphic to $\Z$ as we shall see below, 
which will be identified with the Lie algebra of $k^*$ acting on the hypersuface singularity defined by $W$.
\begin{defn}
Let $W=(a_1,a_2,a_3;h)$ be a regular system of weights of dual type.
Choose a triple $\vec{l}_W:=(\vec{l}_1,\vec{l}_2,\vec{l}_3)$ of elements in $L(A_W)_+$ as follows 
(see also Appendix):
\begin{center}
\begin{tabular}{c||c|c|c}
Type of $W$ & $\vec{l}_1$ & $\vec{l}_2$ & $\vec{l}_3$ \\
\hline
\hline
I & $\vec{X}_1$ & $\vec{X}_2$ & $\vec{X}_3$ \\
\hline
II & $\vec{X}_1+\vec{X}_3$ & $p_1\vec{X}_3$ & $\vec{X}_2$ \\
\hline
III & $\vec{X}_1+\vec{X}_2+\vec{X}_3$ & $p_1\vec{X}_3$ & $p_1\vec{X}_2$ \\
\hline
IV & $\frac{p_3}{p_2}\vec{X}_2+\vec{X}_3$ & $p_1\vec{X}_3$ & $p_1\vec{X}_2$ \\
\hline
V &  $\vec{X}_2+l\vec{X}_3$ & $\vec{X}_1+k\vec{X}_2$ & $\vec{X}_3+m\vec{X}_1$ \\
\end{tabular}
\end{center}
We shall call $\vec{l}_W$ the {\it principal generators} of $W$.
\end{defn}
\begin{lem}\label{lem:omega}
Let $W=(a_1,a_2,a_3;h)$ be a regular system of weights of dual type and $\vec{l}_W=(\vec{l}_1,\vec{l}_2, \vec{l}_3)$ 
be its principal generators.
There exists a unique element $\vec{\omega}_W\in L(A_W)$ such that 
\begin{enumerate}
\item $-\epsilon_W\cdot \vec{\omega}_{W}=\vec{\omega}_{A_W}$
\item $a_i\cdot\vec{\omega}_W=\vec{l}_i$ for $i=1,2,3$.
\end{enumerate}
\end{lem}
\begin{pf}
By direct calculations, one sees that  
\begin{equation}
a_i\cdot\vec{\omega}_{A_W}=-\epsilon_W\cdot \vec{l_i}\quad \text{and}\quad
a_i\cdot\vec{l_j}=a_j\cdot\vec{l_i},\quad \text{for all}\quad i,j=1,2,3.
\end{equation}
Note also that there exists $k_1,k_2,k_3\in\Z$ such that $\sum_{i=1}^3k_ia_i=1$ since ${\rm gcd}(a_1,a_2,a_3)=1$.
(Existence) 
Set $\vec{\omega}_W:=\sum_{i=1}^3k_i\cdot\vec{l_i}$. 
Then 
$$
-\epsilon_W\cdot\vec{\omega}_W=\sum_{i=1}^3k_i(-\epsilon_W\cdot\vec{l_i})=
\sum_{i=1}^3k_i(a_i\cdot\vec{\omega}_{A_W})=(\sum_{i=1}^3k_ia_i)\cdot\vec{\omega}_{A_W}=\vec{\omega}_{A_W}.
$$
Moreover, 
\begin{align*}
a_1\cdot\vec{\omega}_{W}=&a_1k_1\cdot\vec{l_1}+a_1k_2\cdot\vec{l_2}+a_1k_3\cdot\vec{l_3}\\
=&(1-k_2a_2-k_3a_3)\cdot\vec{l_1}+a_1k_2\cdot\vec{l_2}+a_1k_3\cdot\vec{l_3}\\
=&\vec{l}_1+k_2(a_1\cdot\vec{l_2}-a_2\cdot\vec{l_1})+k_3(a_1\cdot\vec{l_3}-a_3\cdot\vec{l_1})\\
=&\vec{l}_1.
\end{align*}
In a similar way, it is shown that $a_2\cdot\vec{\omega}_{W}=\vec{l_2}$ and $a_3\cdot\vec{\omega}_{W}=\vec{l_3}$. 
(Uniqueness)
Suppose that an element $\vec{\omega}'_W\in L(A_W)$ satisfies the two conditions (i) and (ii) above.
Then, 
$$
\vec{\omega}'_W=(\sum_{i=1}^3k_ia_i)\cdot\vec{\omega}'_W=\sum_{i=1}^3k_i(a_i\cdot\vec{\omega}'_W)
=\sum_{i=1}^3k_i\cdot\vec{l_i}=\vec{\omega}_W.
$$
\qed
\end{pf}
\begin{rem}
Let $W=(a_1,a_2,a_3;h)$ be a regular system of weights of genus $0$ and $A_W=(\alpha_1,\dots, \alpha_r)$ be its signature. 
Then one has
\begin{equation}
{\rm deg}(\vec{\omega}_{W})=-\frac{{\rm deg}(\omega_{A_W})}{\epsilon_W}=
-\frac{(r-2)\cdot{\rm deg}(\vec{c})-\sum_{i=1}^r{\rm deg}(\vec{X_i})}{\epsilon_W}
=\frac{h}{a_1a_2a_3}\cdot \alpha,
\end{equation}
where $\alpha:={\rm lcm}(\alpha_1,\dots,\alpha_r)$. 
\end{rem}
This element $\vec{\omega}_W\in L(A_W)$ leads to the following definition: 
\begin{defn}\label{def:subring} 
Let $W$ be a regular system of weights of dual type and $\vec{\omega}_W\in L(A_W)$ be the element given by 
Lemma \ref{lem:omega}.
\begin{enumerate}
\item The subgroup $L_W:=\Z\vec{\omega}_W$ of $L(A_W)$ is called the {\it principal lattice} of $W$, 
which we shall often identify with $\Z$.
\item Define an $L_W$-graded $k$-algebra $R_W$ by 
$$
R_W:=\bigoplus_{d\in\Z_{\ge 0}}R_{W,d},\quad R_{W,d}:=R_{A_W,d\cdot\vec{\omega}_W}.
$$
\end{enumerate}
\end{defn}
Denote by ${\rm gr}\text{-}R_{W}$ the abelian category of finitely generated $L_W$-graded 
$R_{W}$-modules and by ${\rm tor}\text{-}R_{W}$ the full abelian subcategory of 
${\rm gr}\text{-}R_{W}$ whose objects are finite dimensional $L_W$-graded $R_{W}$-modules.
\begin{defn}
Define an algebraic stack $\CC_{W}$ by 
\begin{equation}
\CC_{W}:=\left[{\rm Spec}(R_{W})\backslash\{0\}/{\rm Spec}({k L_W})\right].
\end{equation}
The abelian category ${\rm coh}(\CC_{W})$ of coherent sheaves on the algebraic stack $\CC_{W}$ is equivalent to 
${\rm gr}\text{-}R_{W}/{\rm tor}\text{-}R_{W}$.
\end{defn}
Our main purpose is to compare $\CC_{A_W}$ with $\CC_W$ as algebraic stacks, or equivalently, 
${\rm mod}_{L(A_W)}\text{-}R_{A_W,\lambda}/{\rm tor}_{L(A_W)}\text{-}R_{A_W,\lambda}$ 
with ${\rm gr}\text{-}R_{W}/{\rm tor}\text{-}R_{W}$ as abelian categories, and 
relate them to the hypersuface singularity with $k^*$-action defined by $W$.
\section{Main Theorem}\label{sec:thm}
\begin{thm}\label{thm:main}
Let $W$ be a regular system of weights of dual type.
\begin{enumerate}
\item There exist homogeneous elements $x,y,z\in R_W$ of degree $a_1, a_2, a_3$ with respect to $L_W$
such that $R_W\simeq k[x,y,z]/(f_W)$ for some weighted homogeneous polynomial $f_W\in k[x,y,z]$ of degree $h$
which defines at most an isolated singularity only at the origin.
In other words, $R_W$ is the graded ring of functions on the isolated hypersurface singularity defined by $W$. 
\item There exists an equivalence of abelian categories$:$
\begin{equation}
{\rm mod}_{L(A_W)}\text{-}R_{A_W}\left/{\rm tor}_{L(A_W)}\text{-}R_{A_W}\right.
\simeq {\rm gr}\text{-}R_W\left/{\rm tor}\text{-}R_W\right..
\end{equation}
\end{enumerate}
\end{thm}
\begin{pf}
Proofs for both of two statements are generalizations of those given by Geigle-Lenzing \cite{gl:2} 
for the case $\epsilon_W=1$, which can be applied to the general cases here.
\begin{enumerate}
\item The statement is shown by direct calculations using the classification of regular systems of weights of dual type. 
Generators $x,y,z$ and their unique relation $f_W(x,y,z)$ are listed in Appendix below.
\item One sees that the functor 
\begin{equation}
F:{\rm mod}_{L(A_W)}\text{-}R_{A_W}\to {\rm gr}\text{-}R_W, 
\quad M\mapsto M|_{L_W},
\end{equation}
defined by the natural inclusion $R_W\subset R_{A_W}$ is exact 
and essentially surjective since by Kan extension there exists the right adjoint $G$ of $F$ 
such that $F\circ G (M')\simeq M'$ for any $M'\in {\rm gr}\text{-}R_W$.
Therefore, only one has to to show is that any module $N\in {\rm tor}\text{-}R_W$ is an essential 
image of a module $N'\in{\rm tor}\text{-}R_{A_W}$, however, this also follows from Proposition 1.3 in \cite{gl:1}.
\end{enumerate}
\qed
\end{pf}
\section{Application}\label{sec:cor}
\begin{cor}
Let $W$ be a regular system of weights of dual type with $\epsilon_W<0$.
The triangulated category
$D^{gr}_{Sg}(R_W):=D^b({\rm gr}\text{-}R_W)\left/D^b({\rm grproj}\text{-}R_W)\right.$
has a full exceptional collection $(E_1,\dots,E_{\mu_{W^*}})$.
\end{cor}
\begin{pf}
Orlov shows (Theorem 2.5 in \cite{o:1}) that there exists a semi-orthogonal decomposition
\begin{equation}
D^{gr}_{Sg}(R_W)\simeq \left<R_W/{\mathfrak m}_W(a),\dots,R_W/{\mathfrak m}_W(a+\epsilon_W+1),
D^b{\rm coh}(\CC_W)\right>
\end{equation}
for any $a\in\Z$ where we denote by $(\cdot)$ the autoequivalence induced by 
the grading shift functor on ${\rm gr}\text{-}R_W$ and by 
${\mathfrak m}_W$ the unique graded maximal ideal in $R_W$.
Combining our main theorem with the Proposition 4.1 in \cite{gl:1}, the statement follows.
\qed
\end{pf}
\begin{rem}
There is a isomorphism of functors $(h)\simeq [2]$ 
on the triangulated category $D^{gr}_{Sg}(R_W)$ where $[1]$ is the translation functor.
\end{rem}
\begin{cor}\label{cor:exceptional}
Let $W$ be a regular system of weights of dual type.
The triangulated category $D^{gr}_{Sg}(R_W)$ has a full exceptional collection 
$(E_1,\dots,E_{\mu_{W^*}})$.
\end{cor}
\begin{pf}
Since $W$ is of dual type, it is of ADE type or $\epsilon_W<0$. 
We know that a regular system of type ADE has a full strongly exceptional collection \cite{kst:1},
it is the direct consequence of the previous corollary. 
\qed
\end{pf}
\section{Conjectures}\label{sec:conj}
In this section, we shall assume that $k=\C$.
\begin{defn}Let $W$ be a regular system of weights of dual type.
\begin{enumerate}
\item Set $B_W:=A_{W^*}$ and call it the \it{dual signature} of $W$.  
\item Let $B_W=(\beta_1,\beta_2,\beta_3)$ be the dual signature of $W$. 
Denote by $f_{B_W}$ the defining polynomial of the cusp singularity, namely,   
$f_{B_W}:=x_1^{\beta_1}+x_2^{\beta_2}+x_3^{\beta_3} + x_1x_2x_3$. 
\end{enumerate}
\end{defn}
\begin{rem}
For singularities associated to regular system of weights with $\epsilon_W=-1$, i.e., 
Arnold's 14 exceptional singularities, $B_W$ is called the \it{Gabrielov number}.
\end{rem}
Kontsevich observed that the origin of mirror symmetry is a symmetry between complex geometry 
and symplectic geometry. 
In particular, he conjectures that the algebraically constructed triangulated category in complex side 
should be triangulated equivalent to the geometrically constructed one in symplectic side, and  
that mirror phenomena should be explained from this triangulated equivalence.
Let $f$ be a polynomial which defines an isolated singularity only at the origin $0\in\C^n$. 
Consider a polynomial $f$ as a holomorphic map $f:(\C^n,0)\to (\C,0)$ and choose a 
{\it Morsification} $\{f_{t}\}_{0\le t<1}$ of $f_{}$, i.e., a smooth family of polynomials with 
$f_{0}=f$ and such that critical points of $f_t$ ($t\ne 0$) are isolated and non-degenerate.  
For very small $0<t\ll \delta\ll\epsilon\ll 1$, put 
\begin{equation}
X:=\{ x\in \C^{n}~|~||x||\le\epsilon, |f_t|\le\delta \},
\end{equation}
\begin{equation}
D:=\{y\in\C~|~|y|\le \delta\}.
\end{equation}
Then, after suitable modifications near the boundary, the holomorphic map 
$f_t:X\to D$ defines an {\it exact Lefschetz fibration} and the relative holomorphic form 
$n-1$-form $dx_1\wedge \dots \wedge dx_n/df_t$ defines a {\it relative Maslov map}  
which enables one to define {\it graded Lagrangian submanifolds}.
\begin{thm}[Seidel \cite{se:1} \cite{se:2}]
There exists an $A_\infty$-category ${\rm Fuk}(f)$
called the Fukaya category of the exact Lefschetz fibration such that 
\begin{equation}
D^b{\rm Fuk}(f)\simeq D^b{\rm Fuk}^{\rightarrow}(\{\gamma_\bullet\}),
\end{equation}
as a triangulated category, 
where ${\rm Fuk}^{\rightarrow}(\{\gamma_\bullet\})$ is the directed Fukaya category 
of a distinguished basis of vanishing graded Lagrangian submanifolds
$\{\gamma_\bullet\}:=\{\gamma_1,\dots, \gamma_\mu\}$ in the Milnor fiber $X_y=f_t^{-1}(y)$.
In particular, $D^b{\rm Fuk}(f)$ is independent of all choices and is an invariant of $f$.
\qed
\end{thm}
We expect that the topological mirror symmetry for isolated hypersurface singularities 
can also be "categorified". The homological mirror symmetry principle leads the following conjecture:
\begin{conj}[Homological Mirror Symmetry for Weighted Projective Lines]
Let $W$ be a regular system of weights of dual type.
There should exist a triangulated equivalence
\begin{equation}
D^b{\rm coh}(\CC_W)\simeq D^b{\rm Fuk}(f_{B_W}).
\end{equation}
\end{conj}
\begin{rem}
The above conjecture can be verified at the level of the Grothendieck group. More precicely, 
one has the following isomorphism as a lattice
\begin{equation}
\left(K_0(D^b{\rm coh}(\CC_W)),\chi+{}^t\chi\right)\simeq
\left(H_2(f_{B_W,t}^{-1}(y),\Z),-I\right).
\end{equation}
\end{rem}
We expect that this conjecture follows from the following one since semi-orthogonal decompositions of 
triangulated categories in algebraic geometry side should correspond to operations "blowing  down" 
which are mirror dual to unfoldings of singularities in symplectic geometry side. 
\begin{conj}[HMS for regular systems of weights of dual type]
Let $W$ be a regular system of weights of dual type.
\begin{enumerate}
\item There should exist a triangulated equivalence
\begin{equation}
D^{gr}_{Sg}(R_W)\simeq D^b{\rm Fuk}(f_{W^*}).
\end{equation}
\item $D^b{\rm Fuk}(f_{W^*})$ has the following semi-orthogonal decomposition
\begin{equation}
D^b{\rm Fuk}(f_{W^*})\simeq \left<\L(a),\dots,\L(a+\epsilon_W+1),
D^b{\rm Fuk}(f_{B_W})\right>,
\end{equation}
for any $a\in\Z$, where $\L$ is an object of $D^b{\rm Fuk}(f_{W^*})$ and $(1)$ is the auto-equivalence 
induced by the shift of the gradings of Lagrangian submanifolds by $2/h$ such that $(h)=[2]$ where 
$[1]$ is the translation functor on $D^b{\rm Fuk}(f_{W^*})$.
\end{enumerate}
\end{conj}
\begin{rem}
We can define an abelian group $L_{f_W}$ for $R_W$ called the {\it maximal grading} of $R_W$ in a similar way to Definition \ref{defn:L} 
by introducing the degree vectors for $x_i$ and letting all monomials appearing in $f_W$ have the same degree vectors.
As a result, $R_W$ becomes $L_{f_W}$-graded and hence we can consider the triangulated category $D^{L_{f_W}}_{Sg}(R_W)$.
It is expected that 
if $f_W$ is of the form listed in Appendix (without any conditions for their exponents) 
then $(W,{\rm Spec}(kL_W))$ is topological mirror dual to $(W^*,\{1\})$ with the same $W^*$ in the list.
In this case, we can formulate the homological mirror symmetry conjecture 
as a triangulated equivalence $D^{L_{f_W}}_{Sg}(R_W)\simeq D^b{\rm Fuk}(f_{W^*})$.
The above conjecture is the special case for regular systems of weights with $L_W=\Z$ (see Definition \ref{defn:dual}). 
\end{rem}
This homological mirror symmetry conjecture is one of motivations of our study of the triangulated category 
$D^{gr}_{Sg}(R_W)$, from which we can expect that $D^{gr}_{Sg}(R_W)$ has a full strongly exceptional collection 
(see \cite{t:2}). 
Indeed, we prove this for several important classes of regular systems of weights in \cite{kst:1} (for $W$ corresponding to 
ADE singularities), \cite{kst:2} (for $W$ with $\epsilon_W=-1$ including Arnold's 14 exceptional singularities) 
and \cite{t:3} (for $W$ of Type I, II and III), 
where we also checked the homological mirror symmetry conjecture at the level of Grothendieck group.
\section{Appendix}
The followings are data of weights, dual weights, signatures, generators and relations of $R_W$ and of $R_{W^*}$,
for regular systems of weights of dual type.
\subsection{Type I}
\begin{equation}
W=W^*=(p_2p_3,p_3p_1,p_1p_2;p_1p_2p_3), \quad (p_i,p_j)=1, i=1,2,3.
\end{equation}
\begin{equation}
A_W=A_W^*=(p_1,p_2,p_3).
\end{equation}
\begin{equation}
\vec{l}_W=(\vec{X}_1,\vec{X}_2,\vec{X}_3).
\end{equation}
\begin{equation}
x=X_1,\quad y=X_2,\quad z=X_3.
\end{equation}
\begin{equation}
f_W(x,y,z)=x^{p_1}+y^{p_2}+z^{p_3}.
\end{equation}
\smallskip
\subsection{Type II}
\begin{equation}
W=(p_3,\frac{p_1p_3}{p_2},(p_2-1)p_1;p_1p_3),
\end{equation}
\begin{equation}
W^*=(p_3,p_1p_2,(\frac{p_3}{p_2}-1)p_1;p_1p_3), 
\end{equation}
where $p_2\ne p_3$, $p_2|p_3$, $(p_1,p_3)=1$,
$(p_2-1,p_3)=1$ and $(p_3/p_2-1,p_3)=1$.
\begin{equation}
A_W=(p_1,\frac{p_3}{p_2},(p_2-1)p_1),\quad A_{W^*}=(p_1, p_2, (\frac{p_3}{p_2}-1)p_1).
\end{equation}
\begin{equation}
\vec{l}_W=(\vec{X}_1+\vec{X}_3,p_1\vec{X}_3,\vec{X}_2).
\end{equation}
\begin{equation}
x=X_1X_3,\quad y=X_3^{p_1},\quad z=X_2.
\end{equation}
\begin{equation}
f_W(x,y,z)=x^{p_1}+y^{p_2}+yz^{\frac{p_3}{p_2}}.
\end{equation}
\smallskip
\subsection{Type III}
\begin{equation}
W=W^*=(p_2,p_1q_2,p_1q_3;p_1p_2),
\end{equation}
where $(p_1,p_2)=1$, $p_2+1=(q_2+1)(q_3+1)$ and $(q_2,q_3)=1$.
\begin{equation}
A_W=A_{W^*}=(p_1,p_1q_2,p_1q_3).
\end{equation}
\begin{equation}
\vec{l}_W=(\vec{X}_1+\vec{X}_2+\vec{X}_3,p_1\vec{X}_3,p_1\vec{X}_2).
\end{equation}
\begin{equation}
x=X_1X_2X_3,\quad y=X_3^{p_1},\quad z=X_2^{p_1}.
\end{equation}
\begin{equation}
f_W(x,y,z)=x^{p_1}+y^{q_3+1}z+yz^{q_2+1}.
\end{equation}
\smallskip
\subsection{Type IV}
\begin{equation}
W=(\frac{p_3}{p_1},(p_1-1)\frac{p_3}{p_2},p_2-p_1+1;p_3),
\end{equation}
\begin{equation}
W^*=(p_2,(\frac{p_3}{p_2}-1)p_1,\frac{p_3}{p_1}-\frac{p_3}{p_2}+1;p_3),
\end{equation}
where $p_1\ne p_2\ne p_3$, $p_1|p_2$, $p_2|p_3$,
$(p_1-1,p_2)=1$, $(p_2-p_1+1,p_3)=1$, $(p_3/p_2-1,p_3/p_1)=1$ and $(p_3/p_1-p_3/p_2+1,p_3)=1$.
\begin{equation}
A_W=(\frac{p_3}{p_2}, (p_1-1)\frac{p_3}{p_2}, p_2-p_1+1),\quad 
A_W=(p_1, (\frac{p_3}{p_2}-1)p_1, \frac{p_3}{p_1}-\frac{p_3}{p_2}+1).
\end{equation}
\begin{equation}
\vec{l}_W=\left(\frac{p_3}{p_2}\vec{X}_2+\vec{X}_3,p_1\vec{X}_3,\vec{X}_1+\vec{X}_2\right).
\end{equation}
\begin{equation}
x=X_2^{\frac{p_3}{p_2}}X_3,\quad y=X_3^{p_1},\quad z=X_1X_2.
\end{equation}
\begin{equation}
f_W(x,y,z)=x^{p_1}+xy^{\frac{p_2}{p_1}}+yz^{\frac{p_3}{p_2}}.
\end{equation}
\smallskip
\subsection{Type V}
\begin{equation}
W=(lm-m+1,mk-k+1,kl-l+1;klm+1),
\end{equation}
\begin{equation}
W^*=(lm-l+1,mk-m+1,kl-k+1;klm+1).
\end{equation}
where $(lm-m+1,klm+1)=1$C$(mk-k+1,klm+1)=1$, and $(kl-l+1,klm+1)=1$.
\begin{equation}
A_W=(lm-m+1,mk-k+1,kl-l+1),\quad A_{W^*}=(lm-l+1,mk-m+1,kl-k+1).
\end{equation}
\begin{equation}
\vec{l}_W=(\vec{X}_2+l\vec{X}_3,\vec{X}_1+k\vec{X}_2,\vec{X}_3+m\vec{X}_1).
\end{equation}
\begin{equation}
x=X_2X_3^l,\quad y=X_1X_2^{k},\quad z=X_3X_1^m.
\end{equation}
\begin{equation}
f_W(x,y,z)=zx^{k}+xy^{m}+yz^{l}.
\end{equation}

\end{document}